\documentclass[12pt,a4paper]{article}%
\usepackage[utf8]{inputenc}
\usepackage{hyperref}
\usepackage{amsmath}
\usepackage{amsfonts}
\usepackage{amssymb}
\usepackage{xcolor}
\usepackage{graphicx}%
\providecommand{\U}[1]{\protect\rule{.1in}{.1in}}
\newtheorem{theorem}{Theorem}

\newtheorem{claim}{Claim}

\newtheorem{conjecture}[theorem]{Conjecture}
\newtheorem{corollary}[theorem]{Corollary}

\newtheorem{lemma}[theorem]{Lemma}

\newtheorem{observation}[theorem]{Observation}

\newenvironment{proof}[1][Proof]{\noindent\textbf{#1.} }{\ \hfill \rule{0.5em}{0.5em}\bigskip}
\graphicspath{{Slike/}{../Slike/}}

\begin{document}

\title{Two results on rich-neighbour $5$-colouring of subcubic graphs}
\author{Masaki Kashima$^{1}$, Mirko Petru\v{s}evski$^{2}$, Jelena Sedlar$^{3}$,\\[0.1cm]
Riste \v{S}krekovski$^{4}$\\[0.3cm]
{\small $^{1}$ \textit{Keio University, FST, Yokohama, Japan}}\\[0.1cm]
{\small $^{2}$ \textit{Ss. Cyril and Methodius University, FME, Skopje, North
Macedonia}}\\[0.1cm]
{\small $^{3}$ \textit{University of Split, FGAG, Split, Croatia}}\\[0.1cm]
{\small\parbox{0.95\textwidth}{\centering\textit{
$^{4}$ FMF, University of Ljubljana, FIS, University of Novo Mesto, Slovenia}}}}

\date{}
\maketitle

\begingroup
\renewcommand{\thefootnote}{}
\footnotetext{\textit{E-mail addresses:}
\href{mailto:masaki.kashima10@gmail.com}{masaki.kashima10@gmail.com} (M. Kashima),
\href{mailto:mirko.petrushevski@gmail.com}{mirko.petrushevski@gmail.com} (M.~Petru\v{s}evski),
\href{mailto:jelena.sedlar@gradst.hr}{jelena.sedlar@gradst.hr} (J. Sedlar),
\href{mailto:skrekovski@gmail.com}{skrekovski@gmail.com} (R. \v{S}krekovski).}
\endgroup

\begin{abstract}
A rich-neighbour edge-colouring is a proper edge-colouring in which every
non-isolated edge has an adjacent edge whose neighbours have pairwise
distinct colours. Petru\v{s}evski and \v{S}krekovski~\cite{PS2024} conjectured that every
connected subcubic loopless multigraph distinct from $K_4$ admits such a
colouring with five colours. We prove the conjecture for two classes:
claw-free subcubic loopless multigraphs and $2$-connected subcubic
outerplanar loopless multigraphs. The bound is
sharp in both classes, as witnessed by the diamond $K_4-e$.
\end{abstract}



\section{Introduction}

All graphs considered here are finite loopless multigraphs unless stated otherwise.
Degrees count edge multiplicities, and \emph{subcubic} means maximum degree at most $3$.
A $2$-connected graph is understood to have at least three vertices.
A multigraph is claw-free when its underlying
simple graph has no induced claw. Triangles have three distinct vertices,
and adjacent edges are distinct edges with a common end vertex. We follow \cite{BM2008} for general terminology not defined here. Given a proper edge-colouring $\varphi$ of $G$ with
colour set $\{1,\dots ,k\}$, an edge $e$ is said to be \emph{rich} (under
$\varphi$) if all edges adjacent to $e$ receive pairwise distinct colours; isolated and
pendant edges are automatically rich. Writing $E(u)$ for the edges incident
with $u$ and $E(u,v)$ for all edges joining $u$ and $v$, this is equivalent,
for $e=uv$, to $\varphi(E(u)\setminus E(u,v))\cap
\varphi(E(v)\setminus E(u,v))=\varnothing$. A \emph{rich-neighbour $k$}-\emph{colouring} is a
proper edge-colouring with colour set $\{1,\dots ,k\}$ in which every non-isolated edge has at least one rich neighbour.
The least such $k$ is denoted by $\chi_{\mathrm{rn}}^{\prime}(G)$. A strong edge-colouring makes every edge rich~\cite{FouquetJolivet1983}. Rich-neighbour edge-colouring was introduced as a weakening of strong edge-colouring by Petru\v sevski and
\v Skrekovski \cite{PS2024}, who proved that every subcubic loopless multigraph has a
rich-neighbour $7$-colouring and conjectured the following.

\begin{conjecture}
\label{Con_main}Every connected subcubic graph $G\not =K_{4}$ has a rich-neighbour
$5$-colouring.
\end{conjecture}

The conjecture allows parallel edges, as does Theorem~\ref{Tm_main} below. Chen, Nan and Zhou \cite{CNZ2026} proved the conjectured bound with $5$ replaced by $6$
for two subclasses of simple subcubic graphs, namely the claw-free ones (see Theorem~15 in~\cite{CNZ2026}) and those with
$\mathrm{mad}(G)<36/13$ (see Theorem~12 in~\cite{CNZ2026}), where
$\mathrm{mad}(G)$ denotes the maximum average degree of $G$. They in fact prove slightly more, namely
that in these two graph classes six colours suffice  for a
rich-neighbour colouring under which, in addition, every vertex of degree $2$ is incident
with a rich edge. The conjecture also holds for subcubic graphs whose
vertices of degree~$3$ induce a subgraph of maximum degree at
most~$1$~\cite{CTZ2025}, and bounds linear in the maximum degree
are known for loopless multigraphs~\cite{KPSS2026}. In this paper we settle
Conjecture~\ref{Con_main} for two classes. Our first result covers the whole
claw-free subcubic class, including parallel edges.

\begin{theorem}
\label{Tm_main}Every connected claw-free subcubic multigraph $G\not =K_{4}$ has a
rich-neighbour $5$-colouring.
\end{theorem}

Our second result concerns outerplanar loopless multigraphs.

\begin{theorem}
\label{Tm_2conn}Every $2$-connected subcubic outerplanar multigraph has
a rich-neighbour $5$-colouring.
\end{theorem}

Both bounds are sharp: the diamond $K_4-e$ belongs to both classes and
requires five colours. In fact, each class contains infinitely many graphs
requiring five colours, and neither is covered by~\cite{CTZ2025},
as shown in Section~\ref{Sec_concluding}.

It was observed in~\cite{PS2024} that $K_4$ has to be excluded: in a proper $5$-colouring
of $K_{4}$ some two of the six edges receive the same colour and are then non-adjacent. Each
of the four remaining edges is adjacent to both of them and is therefore not rich, so
neither of those two edges has a rich neighbour.
Throughout the constructions, $\{1,\dots,5\}$ is the palette.
Section~\ref{Sec_clawfree} proves Theorem~\ref{Tm_main} by decomposing a
claw-free graph into bricks and links. Section~\ref{Sec_outerplanar} proves
Theorem~\ref{Tm_2conn} using the outer cycle and the matching of chords.
The two constructions are independent; Section~\ref{Sec_concluding}
discusses their sharpness and the remaining conjecture.

\section{Claw-free subcubic multigraphs}\label{Sec_clawfree}

\subsection{Auxiliary results}\label{Sec_auxiliary}

Let $G$ be a connected claw-free subcubic graph distinct from $K_4$.
If $\Delta(G)\leq2$, then $G$ is a single vertex, a path or a cycle
(possibly of length $2$), and it is easy to see that $G$ has a strong
$5$-edge-colouring.
If $G$ contains three parallel edges, these exhaust the degrees of their
end vertices, so connectivity implies that $G$ consists of these three
edges. Giving them distinct colours is a rich-neighbour $5$-colouring.
Henceforth in this section we assume that $\Delta(G)=3$ and that $G$ has
no triple edge.

A subcubic loopless multigraph is claw-free if and only if every vertex
with three distinct neighbours lies in a triangle. The \emph{bricks} are the subgraphs induced by the vertex sets of the
components of the subgraph formed by all edges that belong to a triangle
or have a parallel mate.
We call their edges \emph{internal}, and an edge with exactly one end vertex
in a brick $B$ is an \emph{outgoing edge} of $B$; thus an edge joining two
bricks is an outgoing edge of both.

\begin{lemma}
\label{Lemma_bricks}Let $G$ be a connected claw-free subcubic loopless
multigraph distinct from $K_4$ and containing no triple edge. Then
\begin{itemize}
\item[(a)] its bricks are pairwise vertex-disjoint induced subgraphs, each
being a simple triangle, a diamond $K_4-e$, two vertices joined by two
parallel edges, or a triangle with one doubled edge;
\item[(b)] every vertex outside the bricks has degree at most $2$;
\item[(c)] every brick vertex has internal degree $2$ or $3$, and hence
at most one outgoing edge; a vertex incident with an outgoing edge has
internal degree exactly $2$;
\item[(d)] every edge whose end vertices belong to the same brick is internal.
\end{itemize}
A simple triangle has at most three outgoing edges, a diamond or a double-edge brick
at most two, and a triangle with a doubled edge at most one.
\end{lemma}

\begin{proof}
For (a), let $C$ be a component of the subgraph formed by all edges
that belong to a triangle or have a parallel mate. If $C$ contains a
parallel pair $uv$, each of $u,v$ has at most one further incident edge.
If they have a common neighbour $w$, these edges form a triangle with a
doubled edge, and no further triangle or parallel pair can meet it.
Otherwise, no further edge at $u$ or $v$ belongs to a triangle or has a
parallel mate. Thus $C$ is a triangle with a doubled edge or a double edge.

If $C$ has no parallel edges, it is a connected union of simple triangles.
Two triangles cannot share just one vertex, since its degree would be at
least $4$. Thus, unless $C$ is a single triangle, it contains two triangles
sharing an edge, forming a diamond. The ends of the shared edge have degree
$3$, and the other two vertices are non-adjacent, since otherwise
connectivity would give $G=K_4$.
No further triangle can meet this diamond without exceeding degree $3$,
so $C$ is the diamond. In all four cases the degree bound and $G\ne K_4$
ensure that $C$ is induced in $G$. This proves (a), including
vertex-disjointness.

For (b), a vertex outside the bricks lies in no triangle and has no
parallel edges; if it had degree $3$, it and its three neighbours would
induce a claw. The internal degrees in the four types in (a) are $2$ or
$3$, so the subcubic degree bound gives (c). Finally, (d) holds because
each brick is an induced subgraph of $G$.
The stated bounds on outgoing edges follow by counting the vertices of
internal degree $2$ in each brick.
\end{proof}

Remove all internal brick edges and call each remaining non-trivial component
a \emph{link}. By Lemma~\ref{Lemma_bricks}\,(b) and (c), every remaining degree
is at most $2$, and each brick vertex has remaining degree at most $1$. So the
links are paths and cycles, and a brick vertex of a link is an end vertex of a
path. If a link contains no vertex of a brick, none of its vertices is incident
with a removed edge. It is therefore a whole component of $G$, so connectivity
implies that $G$ itself is a path or a cycle. Then $G$ has a strong
$5$-edge-colouring: colour its edges greedily in any order, each avoiding the
colours of the at most four edges at distance at most two from it in the line
graph, i.e.\ one or two steps away along $G$. Every edge is then rich, as its
at most two adjacent edges are at distance two from each other.
The edgeless case is immediate. Thus, we may assume that $G$ is neither a path
nor a cycle, and hence each link is a path
with at least one end vertex in a brick; its interior vertices lie outside the bricks and have degree $2$ in $G$,
and an end vertex outside the bricks has degree $1$ in $G$.
We write a link as $P=v_0v_1\cdots v_\ell$, with edges $f_i=v_iv_{i+1}$ for
$0\le i\le \ell-1$, where $\ell$ is its length.

\medskip

Let $B$ be a brick of $G$ and fix colours in $\{1,\dots ,5\}$ on the
outgoing edges of $B$; these are called the \emph{outgoing colours} of $B$.
Say that $B$ is \emph{self-sufficient} if its internal
edges admit a $5$-colouring that is \emph{jointly proper} with the outgoing
colours, i.e.\ proper on the internal and outgoing edges of $B$ together, under which
every internal edge and every outgoing edge of $B$ has a rich neighbour among the
internal edges of $B$ alone.

\begin{lemma}
\label{Lemma_selfsuff} A brick $B$ of $G$ is self-sufficient in each of the following
situations:

\begin{itemize}
\item[(a)] $B$ is a triangle with at most $2$ outgoing edges, with arbitrary outgoing colours;

\item[(b)] $B$ is a diamond with at most one outgoing edge, or a triangle
with one doubled edge, with arbitrary outgoing colour;

\item[(c)] $B$ is a triangle with $3$ outgoing edges whose three outgoing colours take exactly
two distinct values;

\item[(d)] $B$ is a diamond with two outgoing edges or consists of two
parallel edges; the outgoing colours are distinct whenever two are present.
\end{itemize}

\end{lemma}

\begin{proof}
Colours are denoted by Greek letters. By Lemma~\ref{Lemma_bricks}\,(c)
and (d), the edges at a vertex of $B$ are internal edges and at most one
outgoing edge; richness is checked by the criterion from the Introduction.

\medskip

\noindent\textbf{Case (a) or (c):} \emph{$B$ is a triangle $xyz$.}
Choose the labels so that $x$ has no outgoing edge in case (a), and its
outgoing colour occurs once in case (c). The outgoing colours take at
most two values, so colour $xy,xz,yz$ with three distinct colours avoiding
them. This is jointly proper with the outgoing colours. Moreover, $xy$
and $xz$ are rich: their adjacent internal edges have distinct colours
avoiding the outgoing colours, and any outgoing edges at their ends have
distinct colours. Every internal or outgoing edge of $B$ is adjacent to
$xy$ or $xz$.

\medskip

\noindent\textbf{Case (b) or (d):} \emph{$B$ is a triangle $xyz$ with $xy$ doubled, or consists of two
parallel edges joining $x$ and $y$.} Give the internal edges distinct colours avoiding the outgoing colours, which is
jointly proper; there is at most one outgoing colour for a triangle with a
doubled edge and there are at most two for a double edge. In the former,
the internal edges and the outgoing edge have pairwise distinct colours;
in the latter, an internal edge is adjacent to the other internal edge and
to at most two outgoing edges, all of distinct colours by the assumption
in (d). So every internal edge is rich, and every internal or outgoing
edge is adjacent to an internal edge.

\medskip

\noindent\textbf{Case (b) or (d):} \emph{$B$ is a diamond on $u,v,w,x$ with edges
$uv,uw,ux,vw,vx$.} By Lemma~\ref{Lemma_bricks}\,(c), outgoing edges
occur only at $w$ and $x$. If there is a single outgoing edge, we may
assume by symmetry that it is at $w$. Label
the colours $\alpha,\beta,\gamma,\delta,\varepsilon$ so that the outgoing
colours at $w$ and $x$, when present, are $\varepsilon$ and $\delta$,
respectively, which is possible by the assumption in (d), and set
\[
\varphi(uv)=\alpha,\quad\varphi(uw)=\beta,\quad\varphi(ux)=\gamma,\quad
\varphi(vw)=\delta,\quad\varphi(vx)=\varepsilon.
\]
The colours at $u,v,w,x$ lie in $\{\alpha,\beta,\gamma\}$,
$\{\alpha,\delta,\varepsilon\}$, $\{\beta,\delta,\varepsilon\}$ and
$\{\gamma,\varepsilon,\delta\}$, respectively (outgoing colours included),
without repetition, so the colouring is jointly proper. The edges $uv,uw,ux$ are rich, respectively, because
\[
\{\beta,\gamma\}\cap\{\delta,\varepsilon\}=\varnothing,\quad
\{\alpha,\gamma\}\cap\{\delta,\varepsilon\}=\varnothing,\quad
\{\alpha,\beta\}\cap\{\varepsilon,\delta\}=\varnothing.
\]
Every internal or outgoing edge is adjacent to one of them: $uv$ to $uw$,
the other internal edges to $uv$, and the outgoing edges at $w,x$ to
$uw,ux$.
\end{proof}

\medskip

The following classical two-colouring result will supply the outgoing
colours (see~\cite{BM1976}, Lemma 6.1.1, and~\cite{PZ2009}, Claim 4).
The colouring is not required to be proper; we include an Euler-tour proof.

\begin{lemma}
\label{Lemma_polychromatic}
A connected loopless multigraph $J$ admits an edge-colouring with two colours
in which both colours occur at every vertex of degree at least $2$ if and
only if $J$ is not an odd cycle.
\end{lemma}

\begin{proof}
An odd cycle cannot have such a colouring. Conversely, let $J$ be connected
and not an odd cycle. If $J$ has no edges, there is nothing to prove.
Suppose first that all vertices of $J$ have even degrees. If $J$ is a
cycle, it is an even cycle and an alternating colouring works. Otherwise,
choose a vertex $v_0$ of degree at least $4$ and start an Euler tour there.
Colour its edges alternately $1,2$ in their order on the tour.
Every vertex, including $v_0$, occurs internally on the tour, and at
an internal occurrence it sees two consecutive, differently coloured edges.

If $J$ has odd-degree vertices, add a new vertex $v_0$ joined once to each
of them. Since the number of odd-degree vertices is even, the resulting
multigraph is connected and has all degrees even.
Start an Euler tour at $v_0$ and colour its edges alternately $1,2$.
At an even-degree vertex of $J$, any traversal uses two original edges.
At an odd-degree vertex of degree at least $3$, at most one traversal uses
the added edge, so another traversal uses two original edges. In either
case these original edges receive different colours. Deleting the added
edges therefore gives the required colouring of $J$.
\end{proof}

We shall use the following consequence of Lemma~\ref{Lemma_polychromatic}.

\begin{corollary}
\label{Lemma_euler}Let $J$ be a loopless multigraph with $\Delta(J)\le3$. Then $E(J)$ can
be coloured with at most $3$ colours in such a way that every vertex of degree $2$ or $3$
sees exactly two distinct colours on its incident edges.
\end{corollary}

\begin{proof}
We may assume that $J$ is connected. If $J$ is an odd cycle, any proper
$3$-edge-colouring of $J$ has the required property. Otherwise, the
two-colouring supplied by Lemma~\ref{Lemma_polychromatic} has it.
\end{proof}

\subsection{Proof of Theorem \ref{Tm_main}}

\begin{proof}
Let $G$ be a connected claw-free subcubic graph distinct from $K_4$.
The case $\Delta(G)\leq2$ and the triple-edge case were settled in
Subsection~\ref{Sec_auxiliary}, so assume that neither of them occurs.
By the discussion there, every link is a path with at least one end vertex
in a brick, its interior vertices lie outside the bricks and have degree $2$
in $G$, and an end vertex outside the bricks has degree $1$ in $G$.
Notice that no edge of a link has a parallel mate, since such an edge would
be internal.

\begin{claim}
\label{Claim_outgoing}There exists a colouring of the outgoing edges so that
every brick satisfies the requirements of Lemma~\ref{Lemma_selfsuff}.
Moreover, the two outgoing colours of every link of length at least $2$
with both end vertices in bricks are distinct.
\end{claim}

\begin{proof}[Proof of Claim \protect\ref{Claim_outgoing}]
Let $G^{*}$ be obtained from $G$ by contracting each brick $B$ to a vertex
$x_B$, discarding its internal edges, and replacing each link of length
at least $3$ by two edges through a new middle vertex, retaining its
ends. An end in a brick $B$ is represented by $x_B$, and an end outside
the bricks remains a vertex of degree $1$. If both ends lie in the same
brick $B$, the replacement consists of two parallel edges between $x_B$
and the new middle vertex. Links of length $1$ or $2$ are retained after
the contractions.

No length-$1$ link has both ends in the same brick, by
Lemma~\ref{Lemma_bricks}\,(d), so $G^{*}$ is loopless. The degree of $x_B$
is the number of outgoing edges of $B$, at most $3$, and every other
vertex has degree at most $2$. Corollary~\ref{Lemma_euler} therefore
yields a colouring $\psi:E(G^{*})\to\{1,2,3\}$ in which every vertex of
degree $2$ or $3$ sees exactly two distinct colours. Transfer the colour
of each edge at $x_B$ to the corresponding outgoing edge of $B$ in $G$.
For a length-$1$ link joining two bricks, both ends refer to the same
edge of $G^{*}$, so this is well defined.

We check the requirements of Lemma~\ref{Lemma_selfsuff}. If $B$ is a
triangle with $3$ outgoing edges, then $x_B$ has degree $3$ in $G^{*}$, so the
three outgoing colours of $B$ take exactly two distinct values, and case~(c)
applies. If $B$ is a diamond or a double-edge brick with $2$ outgoing edges,
then $x_B$ has degree $2$ in $G^{*}$, so the two outgoing colours of $B$ are
distinct, and case~(d) applies. By Lemma~\ref{Lemma_bricks}, every other
brick is a triangle with at most $2$ outgoing edges, a diamond with at most
one outgoing edge, a triangle with a doubled edge, or a double-edge brick
with at most one outgoing edge, and it is covered by case~(a), (b) or (d)
without any condition on the outgoing colours. Finally, a link of length
at least $2$ with both ends in bricks is represented in $G^{*}$ by two edges
meeting at a vertex of degree $2$. Their colours are distinct, so the two
outgoing colours of the original link are distinct. This establishes the
claim.
\end{proof}

Fix an outgoing-edge colouring from Claim~\ref{Claim_outgoing}. By
Lemma~\ref{Lemma_selfsuff}, extend it to a self-sufficient internal colouring
of each brick. Write $\varphi$ for this partial colouring and its extensions.

\begin{observation}
\label{Obs_rich}Every internal edge and every outgoing edge of a brick has
a neighbour that is an internal edge of the same brick and is rich in $G$,
whatever colours the remaining edges of $G$ receive.
\end{observation}
\noindent Indeed, self-sufficiency of $B$ provides an internal edge $g$ of $B$ with
this property, where richness of $g$ is evaluated on the internal and
outgoing edges of $B$. By Lemma~\ref{Lemma_bricks}\,(c), these are all
edges of $G$ at the end vertices of $g$, so $g$ is rich in $G$.

\begin{claim}
\label{Claim_links}The colouring extends to the interiors of all links so
that it is proper and every non-isolated edge on a link has a rich neighbour.
\end{claim}

\begin{proof}[Proof of Claim \protect\ref{Claim_links}]
Let $P=v_0v_1\cdots v_\ell$ be a link with $v_0$ in a brick and
$f_i=v_iv_{i+1}$. If $\ell=1$, or if $\ell=2$ and both ends are in
bricks, all its edges are already coloured properly by
Claim~\ref{Claim_outgoing} and have rich neighbours by Observation~\ref{Obs_rich}.
If $\ell=2$ and only $v_0$ lies in a brick, give $f_1$ a colour not
appearing on any edge incident with $v_0$. This makes $f_0$ rich, so
$f_1$ has a rich neighbour, while $f_0$ has one by Observation~\ref{Obs_rich}.

Now let $\ell\geq3$. If both ends lie in bricks, the colours of $f_0$
and $f_{\ell-1}$ are already fixed and distinct by
Claim~\ref{Claim_outgoing}. Otherwise, choose the colour of
$f_{\ell-1}$ distinct from that of $f_0$. Let $\alpha,\beta,\gamma$
be the three colours in $\{1,2,3,4,5\}\setminus
\{\varphi(f_0),\varphi(f_{\ell-1})\}$. Since $v_0$ is incident with
exactly two internal brick edges, we may choose the labels so that
$\alpha$ appears on neither of them. For $1\leq i\leq\ell-2$, set
\[
\varphi(f_i)=
\begin{cases}
\alpha & \text{if }i\equiv1\pmod3,\\
\beta & \text{if }i\equiv2\pmod3,\\
\gamma & \text{if }i\equiv0\pmod3.
\end{cases}
\]
The colouring is proper, and $f_0$ is rich by the choice of $\alpha$.
For $1\leq i\leq\ell-2$, the two neighbours $f_{i-1},f_{i+1}$ of
$f_i$ have distinct colours: this follows from the repeating pattern
and the disjoint sets of colours on the end edges and the interior
edges, including $\varphi(f_0)\ne\varphi(f_{\ell-1})$ when $\ell=3$.
Thus every edge except possibly $f_{\ell-1}$ is rich. Each edge
$f_i$ with $i\geq1$ has the rich neighbour $f_{i-1}$, and $f_0$ has
a rich neighbour by Observation~\ref{Obs_rich}. The choices on different links
do not interact, which proves the claim.
\end{proof}

Apply Claim~\ref{Claim_links} to colour all remaining edges. We verify
that the resulting colouring $\varphi$ of $G$ is a rich-neighbour
$5$-colouring. Every edge of $G$ is either internal to a unique brick or
lies on a unique link. At a brick vertex $v$, the edges of $G$ are the
internal edges at $v$ and at most one outgoing edge, by
Lemma~\ref{Lemma_bricks}\,(c), and they are properly coloured by the
definition of self-sufficiency. Every other vertex has degree at most $2$,
its edges lie on one link, and they are properly coloured by
Claim~\ref{Claim_links}. Finally, an internal edge has a rich neighbour by
Observation~\ref{Obs_rich}, and an edge on a link has one by Claim~\ref{Claim_links}.
Hence $\varphi$ is a rich-neighbour $5$-colouring of $G$, which concludes
the proof.
\end{proof}

\section{Two-connected subcubic outerplanar multigraphs}\label{Sec_outerplanar}

In this section parallel edges are allowed. A $2$-connected outerplanar
loopless multigraph has a Hamiltonian outer cycle, and in the subcubic case
the edges outside this cycle form a matching. Some of these edges may be
parallel to cycle edges.
We colour the outer cycle with three colours and the matching with the
remaining two. The chords are coloured by Lemma~\ref{Lemma_tau}, and the
cycle by a rotation of a colouring from Lemma~\ref{Lemma_cycle}; the
non-rich cycle edges produced by the two parts are then kept apart by
Claims~\ref{Claim_spacing} and~\ref{Claim_rotation}. Lemma~\ref{Lemma_cycle}
fails precisely for the orders $4$ and $7$, which are treated explicitly in
Subsection~\ref{Sec_base}.

\subsection{The two colouring lemmas}\label{Sec_prelim}

Fix a $2$-connected subcubic outerplanar loopless multigraph $G$ on
$m\geq3$ vertices. Its underlying simple graph has a unique Hamiltonian
cycle, the boundary of the outer face: a Hamiltonian cycle cannot use an
edge $pq$ joining two nonconsecutive vertices of this boundary, since
$\{p,q\}$ separates the two sides of $pq$. Label the vertices in this cyclic order as $v_0,v_1,\ldots,v_{m-1}$. Choose one copy of
each corresponding edge to form $C=v_0v_1\cdots v_{m-1}$ and draw all
remaining edges inside $C$ without crossings. We call these remaining
edges \emph{chords}, including those parallel to an edge of $C$.
A parallel pair necessarily consists of a cycle edge and a chord: a doubled
edge joining nonconsecutive vertices would give degree at least $4$. Since every vertex already has two cycle edges,
the chords form a matching $M$, so $G=C+M$. Write $e_i=v_iv_{i+1}$,
with indices modulo $m$. The \emph{cyclic distance} between $e_i$ and
$e_j$ is $\min\{|i-j|,m-|i-j|\}$.

\medskip

We seek a proper colouring $\sigma:E(C)\to\{1,2,3\}$ and a colouring
$\tau:M\to\{4,5\}$. Their union is a proper edge-colouring of $G$,
because $M$ is a matching and the palettes are disjoint. Define
\[
\begin{aligned}
N_\sigma&=\{e_i:\sigma(e_{i-1})=\sigma(e_{i+1})\},\\
N_\tau&=\{e_i:e_i\text{ meets two distinct equally coloured chords}\}.
\end{aligned}
\]
Thus $N_\sigma$ is the set of edges not rich in the coloured cycle
$(C,\sigma)$. In the combined colouring, the neighbours of a cycle edge
can repeat a colour only on its two cycle neighbours or on two distinct
incident chords. A chord parallel to the cycle edge is counted only once
among its neighbours, so it creates no additional repetition. Consequently
\begin{equation}
N:=\{e\in E(C):e\text{ is not rich in }G\}
  =N_\sigma\cup N_\tau.\label{For_N}
\end{equation}
A chord not parallel to a cycle edge is not rich under this palette split,
since its four distinct neighbours use only $\{1,2,3\}$. A parallel chord
has three distinct cycle neighbours and may be rich.

For $m\geq5$ with $m\neq7$, our goal is to choose $\sigma$ and $\tau$ so that:

\smallskip
\noindent\textbf{(C1)} No two edges of $N$ have cyclic distance exactly $2$.

\smallskip
\noindent\textbf{(C2)} Every chord is adjacent to at most three edges of $N$.
\smallskip

These conditions will ensure that every cycle edge and every chord has
a rich neighbour. We verify this at the end of the general proof.

\medskip

A chord end vertex is a \emph{foot}. Let the \emph{conflict graph} $H$ be
the multigraph with $V(H)=M$
in which each cycle edge joining feet of two distinct chords joins those
chords. In particular, a cycle edge parallel to a chord contributes no edge
to $H$. We identify the edges of $H$ with their cycle edges in $G$.
Parallel edges form $2$-cycles, and $H$ has no loops by definition. A vertex colouring
of $H$ colours the chords of $G$, with monochromatic edges $N_\tau$.

Every cycle edge lies on exactly one inner face of $G$, and every chord on
exactly two. The weak dual $T$, whose vertices are inner faces and whose edges
represent shared chords, is a tree, since it is connected and each chord
separates the inner faces on its two sides; its vertices include the \emph{digonal} faces
bounded by parallel pairs. Write $C_F$ for the cycle edges
on the boundary of an inner face $F$. A chord $c=pq$ splits $C$ into two
$p$--$q$ paths, the two \emph{sides} of $c$. Since chords do not cross,
every other chord has both feet inside one side. An edge of $H$ not incident
with $c$ has no end in $\{p,q\}$, so it lies on one side of $c$ and joins
two chords of that side.

\begin{observation}
\label{Obs_faces}Every cycle $Z$ of $H$ equals $C_F$ for an inner face $F$
whose boundary alternates between chords and cycle edges. Conversely,
$C_F$ is a cycle of $H$ for every such face with $|C_F|\geq2$.
\end{observation}

\noindent Indeed, let $c=pq$ be a vertex of $Z$. Then $Z-c$ is a path in $H-c$, so
by the above all its vertices are chords on one side $S$ of $c$. Both edges
of $Z$ at $c$ therefore join $c$ to chords on $S$, so they are cycle edges
of $S$. Since only one cycle edge at $p$, and
only one at $q$, lies on $S$, one of them is at $p$ and the other at $q$.
Thus, replacing every vertex of $Z$ by its chord, we obtain a cycle $\hat Z$ of
$G$ alternating between chords and cycle edges, whose vertices are precisely
the feet of its chords. A chord drawn inside $\hat Z$ would have both feet on
$\hat Z$, since $C$ meets the closed region bounded by $\hat Z$ only in
$\hat Z$. This is impossible as $M$ is a matching. Hence $\hat Z$ bounds an inner
face $F$, and $E(Z)=C_F$. Conversely, the boundary of an inner face is a
cycle of $G$, so for $|C_F|\geq2$ its chords are distinct and $C_F$ is a cycle
of $H$. This establishes the observation.

\medskip

Since every cycle edge lies on only one inner face, distinct cycles of $H$
are edge-disjoint. Hence every block of $H$ containing an edge is a
bridge or a cycle, i.e.\ $H$ is a cactus. A digonal face $F$ has
$|C_F|=1$ and gives no cycle of $H$.
Let $\mathcal F$ be the set of inner faces $F$ with alternating boundary
and odd $k=|C_F|\geq3$. By Observation~\ref{Obs_faces}, $F\mapsto C_F$ is a bijection
from $\mathcal F$ onto the set of odd cycles of $H$. The boundary of
$F\in\mathcal F$ alternates between $k$ chords and $k$ cycle edges, so
$F$ has at least six boundary vertices, all of them feet.

Root $T$ at any inner face. For a non-root face, its \emph{parent
chord} is the chord it shares with its parent in $T$. For each
$F\in\mathcal F$, choose $\varepsilon_F\in C_F$ not adjacent to
its parent chord, i.e.\ with no end at a foot of this chord, or any edge
of $C_F$ if $F$ is the root. This is
possible because $|C_F|\geq3$ and only two edges of $C_F$ are adjacent
to the parent chord.

\begin{lemma}
\label{Lemma_tau}The graph $H$ has a vertex $2$-colouring whose monochromatic edges are exactly
$\{\varepsilon_F:F\in\mathcal F\}$. In particular, every even cycle
is properly coloured, every odd cycle has exactly one monochromatic
edge, and each colour class induces a graph of maximum degree at most one.
\end{lemma}

\begin{proof}
Write $D=\{\varepsilon_F:F\in\mathcal F\}$. By Observation~\ref{Obs_faces}, every
cycle of $H$ is $C_F$ for some inner face $F$, and it contains an edge of
$D$ only if $F\in\mathcal F$, because $\varepsilon_{F'}$ lies on no inner
face other than $F'$. Thus $C_F\cap D=\{\varepsilon_F\}$ if $|C_F|$ is
odd and $C_F\cap D=\varnothing$ if $|C_F|$ is even. Hence $H-D$ has no
odd cycle and is bipartite.

Colour $H-D$ properly with two colours. For each $F\in\mathcal F$,
the path $C_F-\varepsilon_F$ has even length and lies in $H-D$, so the
ends of $\varepsilon_F$ receive the same colour. Thus this is a
$2$-colouring of $H$ whose monochromatic edges are exactly $D$.
By the choice of $\varepsilon_F$, $D$ is a matching of $H$. Indeed, as
$M$ is a matching, the ends of $\varepsilon_F$ in $H$ are the two boundary
chords of $F$ next to it, so an end shared by $\varepsilon_F$ and
$\varepsilon_{F'}$, $F\neq F'$, would be a chord on both $F$ and $F'$. But a
chord lies on just two inner faces, a parent and its child in $T$, and the
child does not choose an edge adjacent to that chord. Hence each colour
class induces a graph of maximum degree at most one.
\end{proof}

Fix this colouring of $H$ with colours $4,5$ and let $\tau:M\to\{4,5\}$
be the corresponding chord colouring. Then
\begin{equation}
N_\tau=\{\varepsilon_F:F\in\mathcal F\}.\label{For_Ntau}
\end{equation}
The required cyclic spacing of these designated edges is proved below.

\medskip

\begin{lemma}
\label{Lemma_cycle}If $m\geq3$ and
$m\notin\{4,7\}$, the cycle $C_m$ has a proper edge-colouring with
three colours such that at most two edges are not rich, and if there
are two, their cyclic distance is not $2$.
\end{lemma}

\begin{proof}
Identify the colours with $\mathbb Z_3$. A proper colouring $\sigma$ of
$C_m$ is determined by $\sigma(e_0)$ and the steps
$s_i=\sigma(e_{i+1})-\sigma(e_i)\in\{1,2\}$, which are subject only to the
closing congruence $\sum_{i}s_i\equiv0\pmod 3$. Since
$\sigma(e_{i+1})-\sigma(e_{i-1})=s_{i-1}+s_i$, the edge $e_i$ is rich
precisely when $s_{i-1}=s_i$. Thus the non-rich edges are the positions
where the cyclic step sequence changes its value.

Let $a\in\{0,1,5\}$ with $a\equiv-m\pmod3$. Notice that $a\leq m$, since
$m\equiv1\pmod3$ and $m\notin\{4,7\}$ imply $m\geq10$. Put $s_i=2$ for
$0\leq i<a$ and $s_i=1$ for $a\leq i<m$. Then $\sum_{i}s_i=m+a\equiv0\pmod3$.
If $a=0$, every edge is rich. Otherwise the non-rich edges are $e_0$ and
$e_a$, whose cyclic distance $\min\{a,m-a\}$ equals $1$ for $a=1$ and is
at least $5$ for $a=5$.
\end{proof}

The orders $4$ and $7$ are genuinely exceptional for
Lemma~\ref{Lemma_cycle}, though not for Theorem~\ref{Tm_2conn}. A cyclic sequence taking values in $\{1,2\}$ changes
its value an even number of times. If the steps never change, then
$\sum_is_i\in\{m,2m\}$ forces $3\mid m$. If they change exactly twice, they
form a block of $a$ steps $2$ and a block of $m-a$ steps $1$ with
$a\equiv-m\pmod3$, and the two non-rich edges have cyclic distance
$\min\{a,m-a\}$. For $m=4$ this forces $a=2$, and for $m=7$ it forces
$a\in\{2,5\}$, so in both cases the distance is $2$.

For $m\notin\{4,7\}$, fix a colouring $\sigma$ from Lemma \ref{Lemma_cycle}, relabelling
its colours as $1,2,3$. For $0\leq t<m$, define its rotation by
$\sigma_t(e_i)=\sigma(e_{i-t})$; in particular $\sigma_0=\sigma$.
Each set $N_{\sigma_t}$ has at most two edges and satisfies (C1).

\subsection{The general proof}\label{Sec_2conn}

\begin{proof}[Proof of Theorem \protect\ref{Tm_2conn} for $m\notin\{4,7\}$]
If $m=3$, then $G$ is a triangle with at most one doubled edge.
Giving its at most four edges distinct colours makes every edge rich. Assume
therefore that $m\geq5$ and $m\ne7$, and use the colourings $\tau$
and $\sigma_t$ fixed in Subsection \ref{Sec_prelim}.

\begin{claim}
\label{Claim_spacing}The set $N_\tau$ satisfies
(C1), and every chord is adjacent to at most one edge of $N_\tau$.
\end{claim}

\begin{proof}[Proof of Claim \protect\ref{Claim_spacing}]
We prove the stronger statement that distinct edges of $N_\tau$
have cyclic distance at least $3$. Recall that every cycle edge lies
on exactly one inner face. Let $F\in\mathcal F$, write
$\varepsilon_F=e_i$, and let $e\in\{e_{i-2},e_{i-1}\}$; the case
$e\in\{e_{i+1},e_{i+2}\}$ is symmetric. We first show that $e$ lies on
$F$ or on a proper descendant of $F$ in $T$.

Since the boundary of $F$ alternates between chords and cycle edges,
$v_i$ is a foot of a boundary chord $c=v_iv_j$ of $F$. By the choice
of $\varepsilon_F$, the chord $c$ is not the parent chord of $F$, so
the inner faces on the side of $c$ opposite to $F$ form the subtree of
$T$ rooted at the child of $F$ across $c$. This side contains the
$v_i$--$v_j$ path of $C$ avoiding $e_i$. Notice that $j\neq i+1$, since
otherwise $F$ would be the digonal face bounded by $e_i$ and $c$. If
$j\neq i-1$, the above path has length at least $2$, so it contains
both $e_{i-1}$ and $e_{i-2}$. If $j=i-1$, then $e_{i-1}$ lies on the
digonal face bounded by $e_{i-1}$ and $c$, which is a child of $F$;
moreover, the boundary of $F$ continues after $c$ with a cycle edge at
$v_{i-1}$ distinct from $e_{i-1}$, so $e_{i-2}\in C_F$.

Now let $F,F'\in\mathcal F$ be distinct, and assume to the contrary
that $\varepsilon_F$ and $\varepsilon_{F'}$ have cyclic distance at
most $2$. Since $\varepsilon_{F'}$ lies on $F'\neq F$, the previous
paragraph implies that $F'$ is a proper descendant of $F$. By symmetry,
$F$ is also a proper descendant of $F'$, a contradiction. This proves
the spacing.

Finally, Lemma \ref{Lemma_tau} makes the monochromatic edges a
matching in $H$. Since a vertex of $H$ represents a chord of $G$,
each chord is adjacent to at most one edge of $N_\tau$.
\end{proof}

\begin{claim}
\label{Claim_rotation}There is a rotation $t$
such that $N_{\sigma_t}\cup N_\tau$ satisfies (C1).
\end{claim}

\begin{proof}[Proof of Claim \protect\ref{Claim_rotation}]
Let $q=|\mathcal F|$ be the number of odd cycles of $H$, so
$|N_\tau|=q$ by (\ref{For_Ntau}). We first show that $4q<m$. This is
immediate if $q=0$, so let $q\geq1$. Each face in $\mathcal F$ has at
least six boundary vertices, all of them feet. A foot lies on exactly
two inner faces, namely those on the two sides of its chord. Hence a
foot lies on two faces of $\mathcal F$ only if its chord corresponds
to an edge of $T$ joining two faces of $\mathcal F$. Such edges form a
forest on $\mathcal F$, so there are at most $q-1$ of them, and each
accounts for two feet. Counting distinct feet gives
$m\geq6q-2(q-1)=4q+2$.

Let $X$ be the set of cycle edges at cyclic distance exactly $2$ from
an edge of $N_\tau$, so $|X|\leq2q$. Notice that
$N_{\sigma_t}=\{e_{j+t}:e_j\in N_\sigma\}$. Hence, for $e_j\in
N_\sigma$ and $e_k\in X$, exactly one rotation, namely $t\equiv k-j
\pmod m$, moves $e_j$ to $e_k$. Consequently at most
$|N_\sigma|\,|X|\leq4q<m$ rotations have $N_{\sigma_t}\cap
X\ne\varnothing$, and we choose one of the remaining rotations.
Claim \ref{Claim_spacing} gives (C1) within $N_\tau$, Lemma
\ref{Lemma_cycle} gives (C1) within $N_{\sigma_t}$, and avoiding $X$
rules out distance $2$ between the two sets. Thus their union
satisfies (C1).
\end{proof}

Choose $t$ as in Claim \ref{Claim_rotation} and colour $C$ by
$\sigma_t$ and $M$ by $\tau$. This is a proper colouring with at
most five colours. By Claim \ref{Claim_spacing}, each chord is
adjacent to at most one edge of $N_\tau$, and it is adjacent
to at most $|N_{\sigma_t}|\leq2$ edges of $N_{\sigma_t}$.
Thus $N=N_{\sigma_t}\cup N_\tau$ satisfies (C2) as well as (C1).

By (\ref{For_N}), the cycle edges outside $N$ are rich in $G$.
The two cycle neighbours $e_{i-1},e_{i+1}$ of $e_i$ have cyclic
distance $2$, since $m\geq5$, so (C1) leaves at least one outside
$N$. Every chord not parallel to a cycle edge has four distinct cycle
neighbours, and (C2) leaves at least one outside $N$. A chord parallel to
$e_i$ is adjacent to $e_{i-1}$ and $e_{i+1}$, whose cyclic distance is $2$;
by (C1), at least one of these two edges lies outside $N$ and is rich.
Hence every edge has a rich
neighbour, completing the proof for $m\notin\{4,7\}$.
\end{proof}

\subsection{The exceptional lengths}\label{Sec_base}

\begin{lemma}
\label{Lemma_base}Every $2$-connected subcubic
outerplanar loopless multigraph on $4$ or $7$ vertices admits a rich-neighbour
$5$-colouring.
\end{lemma}

\begin{proof}
Let $m\in\{4,7\}$. One colouring serves both lengths, except for a
single graph on $7$ vertices, up to isomorphism. Colour the cycle edges by
\[
(\sigma(e_0),\ldots,\sigma(e_{m-1}))=
\begin{cases}
(1,2,3,4) & \text{if }m=4,\\
(1,2,3,1,2,3,4) & \text{if }m=7,
\end{cases}
\]
and give every chord colour $5$. Denote this colouring by $\varphi_0$.
Any three consecutive cycle edges receive three distinct colours from
$\{1,2,3,4\}$, so $\varphi_0$ is proper. We need two observations.

First, if a chord is parallel to a cycle edge $e_i$, then both this
chord and $e_i$ are rich, since their neighbours are, respectively,
$e_{i-1},e_i,e_{i+1}$ and $e_{i-1}$, the chord, $e_{i+1}$.

Second, a cycle edge $e_i$ is not rich if and only if $v_i$ and
$v_{i+1}$ are feet of two distinct chords. Indeed, $e_{i-1}$ and
$e_{i+1}$ have distinct colours from $\{1,2,3,4\}$, so the only possible
repetition among the neighbours of $e_i$ is colour $5$ on two distinct
chords at $v_i$ and $v_{i+1}$. Denote by $N_0$ the set of non-rich cycle
edges.

Recall also that, since the chords form a non-crossing matching, for
every chord $uv$ both feet of every other chord lie in the same
component of $C-\{u,v\}$.

Let us show that every cycle edge $e_i$ has a rich neighbour. This is
immediate unless both $e_{i-1}$ and $e_{i+1}$ belong to $N_0$. In that
case the four distinct vertices $v_{i-1},v_i,v_{i+1},v_{i+2}$ are feet.
Their mates are four distinct vertices, and only $m-4\leq3$ vertices
lie outside this set, so some chord joins two of the four vertices. It
is neither $v_{i-1}v_i$ nor $v_{i+1}v_{i+2}$, because $e_{i-1},e_{i+1}\in
N_0$. A chord $v_{i-1}v_{i+1}$ would leave the foot $v_i$ alone in a
component of $C-\{v_{i-1},v_{i+1}\}$, so $v_i$ would have no mate, and
a chord $v_iv_{i+2}$ is excluded in the same way. A chord
$v_{i-1}v_{i+2}$ leaves $v_i,v_{i+1}$ as the only vertices of a
component of $C-\{v_{i-1},v_{i+2}\}$, so they are mates. Thus in every
case $v_iv_{i+1}$ is a chord, and it is a rich neighbour of $e_i$ by
the first observation.

Now consider a chord $uv$. If it is parallel to a cycle edge, then this
edge is a rich neighbour of $uv$ by the first observation. Otherwise
the neighbours of $uv$ are the four cycle edges at $u$ and $v$. Assume
that all of them belong to $N_0$. Then the cycle neighbours of $u$ and
$v$, which are the end vertices of the two components of $C-\{u,v\}$,
are feet of chords distinct from $uv$, with mates in the same
component. Thus no component is a single vertex, which would have to
be its own mate. As the two
components have $m-2\in\{2,5\}$ vertices in total, we obtain $m=7$ and components $x_1x_2$ and $y_1y_2y_3$, where
$x_1,y_3$ are the cycle neighbours of $v$ and $x_2,y_1$ those of $u$,
say. Hence $x_1x_2$ is a chord. If $y_1y_2$ were a chord, then the foot
$y_3$ would have no mate, and similarly $y_2y_3$ is not a chord; hence
$y_1y_3$ is a chord. Therefore $M=\{uv,x_1x_2,y_1y_3\}$.

We conclude that $\varphi_0$ is a rich-neighbour $5$-colouring of $G$,
unless $m=7$ and $G$ is the graph just described. Since $\varphi_0$
was defined for an arbitrary labelling of $C$, in the exceptional case
we may choose the labelling so that
$u=v_2$, $v=v_6$, $x_1x_2=v_0v_1$ and $y_1y_3=v_3v_5$, i.e.
\[
M=\{v_0v_1,v_2v_6,v_3v_5\}.
\]
Recolour the chord $v_3v_5$ with colour $4$. The colouring stays
proper, since colour $4$ appears on the cycle only at $e_6=v_6v_0$,
which is not adjacent to $v_3v_5$. Now the edges $e_0,e_2,e_3,e_4$ and
$v_0v_1$ are rich, as their neighbours have the colours $\{4,5,2\}$,
$\{2,5,1,4\}$, $\{3,4,2\}$, $\{1,4,3\}$ and $\{4,1,2\}$, respectively.
The edge $e_0$ is adjacent to $v_0v_1$, the edges $e_2,e_3,e_4$ are
consecutive, and every other edge is adjacent to one of $e_0,e_2,e_3,e_4$,
which proves the lemma.
\end{proof}

Lemma \ref{Lemma_base} proves Theorem \ref{Tm_2conn} for the two
remaining lengths. Together with Subsection \ref{Sec_2conn}, this
completes the proof of the theorem.

\section{Concluding remarks}
\label{Sec_concluding}

Theorems~\ref{Tm_main} and~\ref{Tm_2conn} confirm Conjecture~\ref{Con_main} for claw-free subcubic loopless multigraphs and $2$-connected subcubic outerplanar loopless multigraphs. The first result covers the class
of claw-free subcubic loopless multigraphs and lowers the six-colour bound of~\cite{CNZ2026} to five, which is best
possible as shown below, although our colourings need not have the additional
property of~\cite{CNZ2026} concerning vertices of degree~$2$. The multigraphs
in this class are exactly those in which every vertex with three distinct
neighbours lies in a triangle. The proof exploits the rigid brick types
supplied by Lemma~\ref{Lemma_bricks}, while all remaining vertices have
degree at most $2$. In the simple case the bricks are precisely triangles
and diamonds.

Five colours are necessary even for the simple diamond $K_4-e$. To see this,
write its edges as $ab,ac,ad,bc,bd$ and suppose that it has a rich-neighbour
$4$-colouring. The four edges other than $ab$ use at most three colours,
so $ab$ is not rich. For $ac$ to have a rich neighbour, one of $ad,bc$ must be
rich, forcing $\varphi(ac)\ne\varphi(bd)$. Similarly, the condition at $ad$
forces $\varphi(ad)\ne\varphi(bc)$. Together with properness these inequalities
make $ac,ad,bc,bd$ pairwise differently coloured, a contradiction.
The diamond is also $2$-connected and outerplanar. Thus the bound is sharp in both classes, already for simple graphs.

Both bounds are attained by infinite families of simple graphs. In a
simple graph, each rich edge whose end vertices have degree~$3$ has four
distinct neighbours, and hence together with them uses five
colours~(cf.~\cite{CNZ2026}). Thus every connected simple cubic claw-free
graph other than $K_4$ requires five colours. Such graphs of arbitrarily
large order are
obtained by replacing each vertex of a connected simple cubic graph by a
triangle and attaching its three incident edges to distinct triangle
vertices. Theorem~\ref{Tm_main} therefore gives
$\chi_{\mathrm{rn}}^{\prime}(G)=5$ for every graph in this family.

For Theorem~\ref{Tm_2conn}, consider the ladders $P_k\square K_2$ for
$k\geq5$, where $P_k$ is the path on $k$ vertices. These graphs are simple,
$2$-connected, subcubic and outerplanar. The four edges adjacent to the
third rung all join vertices of degree $3$. By the above, none of these
four edges can be rich in a proper colouring with at most four colours, so the third
rung would have no rich neighbour. Thus
$\chi_{\mathrm{rn}}^{\prime}(P_k\square K_2)=5$ for every $k\geq5$.

For the outerplanar result away from the exceptional lengths, the mechanism is the split
of the palette into three colours for the outer cycle and two for the chords.
This split is available whenever a graph is a Hamiltonian cycle together
with a matching, and it makes the failure of richness fall
apart into two sets. The first is $N_\sigma$, governed by the
closing congruence $\sum_{i}s_{i}\equiv 0$ $(\mathrm{mod}$ $3)$ along the
outer cycle; the second is $N_\tau$, the set of monochromatic edges of $H$, with exactly
one chosen in each odd cycle and none elsewhere. What makes the two reconcilable
is that the second set contains one edge per odd cycle of $H$,
chosen using the corresponding faces along a rooted weak dual tree, and that the first set consists of
at most two edges whose positions can be rotated together around the cycle.

Conjecture \ref{Con_main} itself remains open. It would be worthwhile to
consider whether the conflict graph of the present paper, whose cycles
correspond to non-digonal inner faces with alternating boundaries, has an analogue for
the graphs of maximum
average degree less than $36/13$ treated in \cite{CNZ2026}, since
what the outerplanar structure really provides here is the weak dual tree
along which the designation is made.

\bigskip\noindent\textbf{Acknowledgments.}~~The authors acknowledge the partial support by
JSPS KAKENHI JP25KJ2077,
by the Slovenian Research and Innovation Agency program P1--0383, by Project KK.01.1.1.02.0027 co-financed by the European
Regional Development Fund, by the NextGeneration EU foundation via IP--UNIST--17 (GEORAZ),
and by the bilateral Slovenian--Croatian project BI--HR/25--27--004.

\bigskip\noindent\textbf{Declaration of AI use.}~~The authors used Claude.ai and
ChatGPT to assist in the preparation of this manuscript.
The authors take full responsibility for the content of the paper,
including all mathematical statements and proofs.

\end{document}